\documentclass[12pt, a4paper]{amsart}
\usepackage[utf8]{inputenc}
\usepackage{amsfonts,amssymb,amsxtra,amsthm,amsmath,amscd,mathrsfs,cite}
\usepackage[all]{xy}%mj-eucal,eucal 可修改 LaTeX 的数学字体命令 \mathcal。当加载该宏包后，使用 \mathcal 命令，调出的是欧拉书写体，而不是通常的计算机现代书写体。
\usepackage{tikz}

\usepackage[normalem]{ulem}
\usepackage{soul}
\usepackage{color}

\setstcolor{red}

\makeatletter
\@namedef{subjclassname@2010}{%
  \textup{2010} Mathematics Subject Classification}
\makeatother

\def \D {{\mathbb D}}

\def \R {{\mathbb R}}

\def\le{\leqslant}
\def\ge{\geqslant}

\theoremstyle{plain}
\newtheorem{theorem}{Theorem}[section]
\newtheorem{proposition}{Proposition}[section]
\newtheorem{lemma}[proposition]{Lemma}

\theoremstyle{remark}

\numberwithin{equation}{section}

\addtocounter{footnote}{1}

\begin{document}
\vglue -5mm

\title[Three-dimensional exponential sums under constant perturbation]
{Three-dimensional exponential sums under constant perturbation}
\author{Jiamin Li~\& Jing Ma}

\address{%
		Jiamin Li
		\\
		School of Mathematics
		\\
		Jilin University
		\\
		Changchun 130012
		\\
		P. R. China
	}
	\email{lijiamin\_math@163.com}
	
	\address{%
		Jing Ma
		\\
		School of Mathematics
		\\
		Jilin University
		\\
		Changchun 130012
		\\
		P. R. China
	}
	\email{jma@jlu.edu.cn}
	
	\date{\today}

\begin{abstract}
By generalizing Bombieri and Iwaniec's double large sieve inequality,
we obtain an estimation on  three-dimensional exponential sums with constant perturbation.
As an application of our estimation on the three-dimensional exponential sums, we get that for any $\epsilon>0$,
$$
\sum\limits_{n\le x}\Lambda\bigg(\bigg[\dfrac{x}{n}\bigg]\bigg)=x\sum\limits_{d\ge 1}\dfrac{\Lambda(d)}{d(d+1)}+O_{\epsilon}\big(x^{\frac{17}{36}+\epsilon}\big)
$$
as $x\rightarrow\infty$.
This improves a recent result of Liu-Wu-Yang, which requires ${9}/{19}$ in place of ${17}/{36}$.
\end{abstract}

\keywords{Exponential sum, Constant perturbation, Double large sieve inequality}
	
\maketitle
	
\section{Introduction}
	The exponential sums of type
	\begin{equation}\label{ShighD}
	\sum\limits_{n_{1}}\cdots \sum\limits_{n_{j}}e(xn_{1}^{\alpha_{1}}\cdots n_{j}^{\alpha_{j}}),
	\end{equation}
	where $n_{1},\dots,n_{j}$ ranges over integers from an interval, have been widely studied and applied in analytic number theory.
In this paper, we consider the following three-dimensional exponential sums
	\begin{equation}\label{Sdelta}
	S_{\delta} (H,M,N):=\sum\limits_{h\sim H}\sum\limits_{m\sim M}\sum\limits_{n\sim N}a(h,m)b(n)e\bigg(X\dfrac{M^{\beta}N^{\gamma}}{H^{\alpha}}\dfrac{h^{\alpha}}{m^{\beta}n^{\gamma}+\delta}\bigg),
	\end{equation}
where $a\sim A$ indicates $A<a\le 2A$, $e(t)=e^{2\pi it}$,
$H, M, N$ are positive integers,
$X>1$ is a real number,
$a(h,m)$ and $b(n)$ are complex numbers such that $|a(h,m)|, |b(n)|\le 1$,
and $\alpha, \beta, \gamma\in \mathbb{R}$ are fixed positive numbers.
Note that $\delta>0$ is the constant perturbation we are concerned with.
	
Review the history, the aforementioned exponential sum \eqref{Sdelta} with $\delta =0$
has been widely studied and applied.
Applying Bombieri and  Iwaniec's double large sieve inequality \cite{BI86},  
Fouvry and Iwaniec \cite[Theorem 3]{FI89} proved that	for $\alpha\beta(\beta+1)\gamma \neq 0$,
	\begin{equation}\label{S01}	
S_{0}(H,M,N)
\ll_{\epsilon}(HMN)^{1+\epsilon}
\bigg(
\bigg(\dfrac{X}{HMN^{2}}\bigg)^{\frac{1}{4}}      +\dfrac{1}{N^{\frac{3}{10}}}     +\dfrac{1}{(HM)^\frac{1}{4}}      +\dfrac{N^{\frac{1}{10}}}{X^{\frac{1}{4}}}
\bigg).
	\end{equation}
Latter, Sargos and Wu \cite[Theorem 7]{SW} proved that for $\alpha(\alpha-1)(\alpha-2)\beta\gamma \neq 0$,
$$
S_{0}(H,M,N) \hskip-0.3mm
\ll_{\epsilon}(HMN)^{1+\epsilon} \hskip-0.3mm
\bigg(\hskip-0.3mm
\bigg(\dfrac{X^4}{H^4M^4N^{11}}\bigg)^{\frac{1}{26}}  \hskip-1mm    \hskip-0.3mm+\hskip-0.3mm   \bigg(\dfrac{X}{HMN^{2}}\bigg)^{\frac{1}{4}}    \hskip-1mm
        +\dfrac{1}{N^\frac{7}{18}}    \hskip-0.3mm+\hskip-0.3mm    \dfrac{1}{(HM)^\frac{1}{4}}    \hskip-0.3mm+\hskip-0.3mm    \dfrac{1}{X^\frac{1}{2}}
\bigg).
$$
Robert and Sargos \cite[Theorem 1]{RS06} proved that	
	\begin{equation}\label{S02}	
S_{0}(H,M,N)
\ll_{\epsilon}(HMN)^{1+\epsilon}
\bigg(       \bigg(\dfrac{X}{HMN^{2}}\bigg)^{\frac{1}{4}}
            +\dfrac{1}{(HM)^{\frac{1}{4}}}
            +\dfrac{1}{N^{\frac{1}{2}}}
            +\dfrac{1}{X^{\frac{1}{2}}}              \bigg).
	\end{equation}

In \cite{LWY1}, Liu, Wu and Yang studied the exponential sums of type \eqref{Sdelta} with constant perturbation $\delta$
and used it to prove an asymptotic formula on a sum involving the Mangoldt function.
In Proposition \ref{prop2.2}, we generalize  Bombieri and Iwaniec's double large sieve inequality
 by generalizing one of the constant sequences to a sequence of functions.
This generalization, combining with the three lemmas we get on Diophantine problems,
finally enable us to prove the following theorem.

\begin{theorem}\label{thm_3DES}
	Let $H, M, N$ be positive integers, $X>1$ be a real number, $a(h,m)$, $b(n)$ be complex numbers such that $|a(h,m)|, |b(n)|\le 1$,
and $\alpha, \beta, \gamma> 0$ be fixed.
Let $\delta>0$ be a   constant perturbation.
If $X\le (8\delta)^{-1}KM^{\beta}N^{\gamma}$ for some $K\ge 1$, % and $M^{\beta}\le \delta N^{\gamma}$,
then
	\begin{equation}\label{SdeltaR}
	S_{\delta}(H,M,N)	\ll_{\epsilon}(HMN)^{1+\epsilon}
 \bigg(
 \bigg(\dfrac{KX}{HMN^{2}}\bigg)^{\frac{1}{4}}      +\Big( \dfrac{K^2}{HM}\Big)^{\frac{1}{4}}    +\Big(\dfrac{K}{N} \Big)^{\frac{1}{2}}   +\dfrac{K}{X^{\frac{1}{2}}}
 \bigg).
	\end{equation}
\end{theorem}
	
Note that \cite[Theorem 1]{RS06} implies \eqref{SdeltaR} for $\delta=0$.

\vskip1mm
%%%%%%%%%%%%%%%%%%%%%%%%%%%%%%%%%%%%%%%%%%%%%%%%%%%%%%%%%%%%%%%%%%%%%%%%%%%%%%%%%%%%%%%%%%%55
As an application of Theorem \ref{thm_3DES}, we consider the asymptotic behaviour of the following sum
$$
S_{\Lambda}(x):=\sum\limits_{n\le x}   \Lambda \bigg(\bigg[\dfrac{x}{n}\bigg]\bigg)
$$
as $x\rightarrow \infty$, where $\Lambda$ is the Mangoldt function.
This type of sum is a generalization of Dirichlet divisor problem,
 and the study on it is initiated in Bordell\`{e}s-Dai-Heyman-Pan-Shparlinski \cite{BDHPS19}.
Many authors studied this type of sum for different arithmetic functions, see \cite{MaWuzhao,LM1,ZW,Wu20,Wu2019,Zh20,MS21,MS22,LWY2,St21}.
For the Mangoldt function,   Bordell\`{e}s-Dai-Heyman-Pan-Shparlinski's result \cite{BDHPS19} implies
\begin{equation}\label{Mangldot-formula}
	S_{\Lambda}(x)=x\sum\limits_{d\ge 1}\dfrac{\Lambda(d)}{d(d+1)}+O_{\epsilon}(x^{2/3+\epsilon}).
\end{equation}
%improved Bordell\`{e}s-Dai-Heyman-Pan-Shparlinski's work and implied that  the exponent in \eqref{Mangldot-formula} can be
Wu \cite{Wu20} reduced the exponent in   \eqref{Mangldot-formula} to $1/2$.
With the help of the Vaughan identity and the method of one-dimensional exponential sum,
Ma and Wu \cite{MW20}  broke the $1/2$-barrier and reduced the exponent in \eqref{Mangldot-formula} to $ 35/71$.
Latter,  using a result of Baker on two-dimensional exponential sums,
Bordell\`{e}s  \cite{B20} sharpened the exponent in \eqref{Mangldot-formula} to $ {97}/{203}$.
Recently, using the method of multiple exponential sums,
Liu-Wu-Yang \cite{LWY1} revised  the exponent in \eqref{Mangldot-formula} to ${9}/{19}$.
Applying our  estimation on the three-dimensional exponential sums with constant perturbation  we get in Theorem \ref{thm_3DES},
we are able to   revise   the exponent in \eqref{Mangldot-formula} to ${17}/{36}$.

\begin{theorem}\label{thm_lambda}
	For any $\epsilon>0$, we have
	$$\sum\limits_{n\le x}\Lambda\bigg(\bigg[\dfrac{x}{n}\bigg]\bigg)=x\sum\limits_{d\ge 1}\dfrac{\Lambda(d)}{d(d+1)}+O_{\epsilon}\big(x^{\frac{17}{36}+\epsilon}\big)\qquad (x\rightarrow\infty).$$
\end{theorem}

\section{A generalization of Bombieri and Iwaniec's double large sieve inequality}

%In this section, we will generalize  Bombieri and Iwaniec's double large sieve inequality.

The exponential sums \eqref{ShighD} can be regraded as a special case of bilinear forms
$$
\mathscr{B}(\mathscr{X},\mathscr{Y}):=
    \sum\limits_{x_r\in \mathscr{X}} \sum\limits_{y_s\in\mathscr{Y}}
a(r)b(s)
e(x_ry_s),
$$
where $\mathscr{X}$, $\mathscr{Y}$ are two set of real numbers with $|x_r|\le X$, $|y_s|\le Y$,
$a(r)$ and $b(s)$ are complex numbers for $x_r\in \mathscr{X}$, $y_s\in \mathscr{Y}$.

To estimate the exponential sums with   perturbation,
 we need to consider the following more general bilinear forms	
$$
\mathscr{B}(a,b;\mathscr{X},\mathscr{Y}):=
%\underset{|\varphi|\le X\ }{\sum\limits_{\varphi\in \mathscr{X}}}\underset{|y|\le Y}{\sum\limits_{y\in\mathscr{Y}}}   a(\varphi)b(y)e(\varphi(y)y),
{\sum\limits_{\varphi\in \mathscr{X}}}   {\sum\limits_{y\in\mathscr{Y}}}   a(\varphi)b(y)e(\varphi(y)y),
$$
where
$a(\varphi), b(y)\in \mathbb{C}$ for $\varphi \in \mathscr{X}$ and $y\in \mathscr{Y}$ respectively;
$\mathscr{Y}$ is a  set of real numbers with $|y|\le Y$ for each $y\in \mathscr{Y}$,
$\mathscr{X}$ is a set of real-valued functions $\varphi\in C([-Y,Y],\R)$ with $|\varphi|\le X$.
Here we have used the notation
$$
|\varphi-\varphi^{*}|:=\sup\limits_{y_{1},y_{2}\in [-Y, Y] }|\varphi(y_{1})-\varphi^{*}(y_{2})|,
$$
in particular, $|\varphi-c|:=\sup\limits_{y\in [-Y, Y]}|\varphi(y)-c|$
for constant $c\in\R$ and $|\varphi|:=\sup\limits_{y\in [-Y, Y]}|\varphi(y)|$.

The following lemma is a direct corollary of \cite[Lemma 2.3]{BI86}.
\begin{lemma}\label{lem2.1}
	Let $\mathscr{Y}$  and $b(y)$ be as above.
 Then we have	
$$
\int_{-T}^{T}  \bigg|  \sum\limits_{y\in\mathscr{Y}}b(y)e(yt)  \bigg|^{2}dt\le (2T+\eta^{-1})\underset{|y-y^*|\le\eta}{\sum\limits_{y\in \mathscr{Y}}\sum\limits_{y^*\in\mathscr{Y}}}|b(y)b(y^*)|
$$
for any positive numbers  $\eta$ and $T$.
\end{lemma}

Using this lemma, we will generalize Bombieri and Iwaniec's double large sieve inequality
and get the following proposition,
which will be applied to prove Theorem \ref{thm_3DES}.

\begin{proposition}\label{prop2.2}
	Let $\mathscr{X}$, $\mathscr{Y}$, $X$, $Y$, $a(\varphi)$,  $b(y)$ be as above, and $K\ge 1$.
	If
	\begin{equation}\label{condition}
	\sup\limits_{y_{1},y_{2}\in [-Y,Y]}|\varphi(y_{1})-\varphi(y_{2})|<\dfrac{K}{4Y}
	\end{equation}
for any $\varphi\in\mathscr{X}$,	
then
	\begin{equation}\label{newBI}
	|\mathscr{B}(a,b;\mathscr{X},\mathscr{Y})|^{2}\ll (1+K  X  Y)\mathscr{B}(b;\mathscr{X})\mathscr{B}(a;\mathscr{Y}),
	\end{equation}
where
$$
\mathscr{B}(b;\mathscr{X}):=
\underset{|y-y^{*}|\le X^{-1}}{\sum\limits_{y\in \mathscr{Y}}\sum\limits_{y^{*}\in\mathscr{Y}}}|b(y)b(y^{*})|,
\qquad	
\mathscr{B}(a;\mathscr{Y}):=
\underset{|\varphi-\varphi^{*}|\le K Y^{-1}}{\sum\limits_{\varphi\in \mathscr{X}}\sum\limits_{\varphi^{*}\in\mathscr{X}}}|a(\varphi)a(\varphi^{*})|.
$$
\end{proposition}

\begin{proof}
Put $\epsilon:=(4Y)^{-1}$, $T:=X+\epsilon$, and
$$
\omega(y):=\dfrac{\pi y}{\sin(2\pi\epsilon y)},\qquad
b^{*}(y):=b(y)\omega(y)
$$
for $y\in\mathscr{Y}$.	
Then
	\begin{equation}\label{wy}
	|\omega(y)|\le \pi Y
	\end{equation}	
for  $|y|\le Y$.
For $|t|\le T$ and $\varphi\in \mathscr{X}$,
if there is a $y\in \mathscr{Y}$ such that $|\varphi(y)-t|\le \epsilon$,
then \eqref{condition} implies $|\varphi-t|\le K/(2Y)$.
Thus, using
\begin{equation*}\label{etydt}
	e(\varphi(y)y)=\dfrac{\pi y}{\sin(2\pi\epsilon y)}\int_{\varphi(y)-\epsilon}^{\varphi(y)+\epsilon} e(ty)dt,
\end{equation*} %with proper consideration of the upper and lower bounds of the integration after interchanging the order of summation and integration, we proceed as follows
we get	
\begin{equation}\label{babxy}
\begin{aligned}
|\mathscr{B}(a,b;\mathscr{X},\mathscr{Y})|
&= \Bigg|{\sum\limits_{\varphi\in \mathscr{X}}}a(\varphi)  {\sum\limits_{y\in\mathscr{Y}}}   b^*(y)
   \int_{\varphi(y)-\epsilon}^{\varphi(y)+\epsilon}   e(ty) dt  \Bigg|
	\\
&= \Bigg|  \int_{-T}^{T}   \underset{\exists\;  y\in \mathscr{Y} \textrm{ such that }   |\varphi(y)-t|\le \epsilon}{\sum\limits_{\varphi\in \mathscr{X}}}   a(\varphi)
                           \underset{|\varphi(y)-t|\le \epsilon}{\sum\limits_{y\in\mathscr{Y} }}   b^*(y)  e(ty) dt  \Bigg|
	\\
&\le \int_{-T}^{T}   \underset{  |\varphi-t|\le    \frac{K}{2Y} }{\sum\limits_{\varphi\in \mathscr{X}}}  | a(\varphi) |
                         \Bigg|  \underset{|\varphi(y)-t|\le \epsilon}{\sum\limits_{y\in\mathscr{Y} }}   b^*(y)  e(ty) \Bigg| dt
	\\
&\le \int_{-T}^{T}   \underset{  |\varphi-t|\le    \frac{K}{2Y} }{\sum\limits_{\varphi\in \mathscr{X}}}  | a(\varphi) |
                         \max_{\mathscr{Y'} \subseteq \mathscr{Y}} \Bigg|   \sum\limits_{y\in\mathscr{Y'}}   b^*(y)  e(ty) \Bigg| dt
	\\
&\le \Bigg( \hskip-1mm   \int_{-T}^{T}   \Big(
         \underset{ |\varphi-t|\le                 \frac{K}{2Y} }{\sum\limits_{\varphi\in \mathscr{X}}}  |a(\varphi)|\Big)^{2} dt\Bigg)^{\frac{1}{2}}
     \Bigg(\int_{-T}^{T}    \Big|     \sum\limits_{y\in\mathscr{Y}_0}    b^{*}(y)e(ty)\Big|^{2}dt\Bigg)^{\frac{1}{2}}
	\end{aligned}
	\end{equation}
by Cauchy's inequality,
where $\mathscr{Y'}$ runs over all the subsets of $\mathscr{Y}$
and $\mathscr{Y}_0$ is a subset of $\mathscr{Y}$ such that the absolute value in the second factor is maximal.
For the first factor, we have
\begin{equation*}
\begin{aligned}
\int_{-T}^{T}\bigg(    \underset{ |\varphi-t|\le      \frac{K}{2Y}}{\sum\limits_{\varphi\in \mathscr{X}}}  |a(\varphi)|   \bigg)^{2}dt
&\le
\int_{-T}^{T}          {\underset{  |\varphi- \varphi^*|\le \frac{K}{Y}, |\varphi-t|\le      \frac{K}{2Y}  } {\sum\limits_{\varphi\in \mathscr{X}}  \sum\limits_{\varphi^*\in \mathscr{X}}}}
|a(\varphi) a(\varphi^*)|      dt
%\\
%&\le
%        \underset{ |\varphi-t|\le 2\epsilon}{\underset{|\varphi|,|\varphi^*|\le X, |\varphi- \varphi^*|\le 4\epsilon}{\sum\limits_{\varphi\in \mathscr{X}}  \sum\limits_{\varphi^*\in \mathscr{X}}}}
%|a(\varphi) a(\varphi^*)|
%\int_{-T}^{T}  dt
\\
&\le 4K \epsilon
\underset{  |\varphi-\varphi^{*}|\le      \frac{K}{Y}  }{\sum\limits_{\varphi\in \mathscr{X}}\sum\limits_{\varphi^{*}\in\mathscr{X}}}|a(\varphi)a(\varphi^{*})|.
\end{aligned}
\end{equation*}
For the second factor,  using Lemma \ref{lem2.1} and \eqref{wy} we get
\begin{align*}
 \int_{-T}^{T}      \Big|     \sum\limits_{y\in\mathscr{Y}_0}
b^{*}(y)e(ty)\Big|^{2}dt
 \le (2T+\eta^{-1})(\pi Y)^{2}\underset{|y-y^{*}|\le\eta}{\sum\limits_{y\in \mathscr{Y}}\sum\limits_{y^{*}\in\mathscr{Y}}}|b(y)b(y^{*})|
\end{align*}
for any $\eta>0$.
Finally,
taking $\eta=X^{-1}$,  \eqref{babxy} implies \eqref{newBI}.
\end{proof}

\section{Proof of Theorem \ref{thm_3DES}}
To prove Theorem \ref{thm_3DES}, we will use Proposition \ref{prop2.2} and  the following  lemmas on Diophantine problems.

\subsection{Some diophantine problems}

The first one is a direct corollary of   \cite[Lemma 1]{FI89}.

\begin{lemma}\label{lem B1}
	Let $\alpha, \beta\in \mathbb{R}^{*}$ be fixed. For $H$, $M\ge 1$ and $X>1$, define
$$
\mathscr{B}_{1}:
=\Big| \Big\{(h_{1},h_{2},m_{1},m_{2})
\ \Big|\ h_{i}\sim H, m_{i} \sim M, i=1, 2,
\Big| \dfrac{h_{1}^{\alpha}m_{1}^{\beta}}{H^{\alpha}M^{\beta}}-\dfrac{h_{2}^{\alpha}m_{2}^{\beta}}{H^{\alpha}M^{\beta}}\Big|\le\dfrac{1}{X}\Big\}\Big|.
$$
Then
$$
\mathscr{B}_{1}\ll_{\epsilon} (HM)^{2+\epsilon}\bigg(\dfrac{1}{HM}+\dfrac{1}{X} \bigg).
$$
\end{lemma}

The second one is  \cite[Theorem 2]{RS06}
which will be used in the proof of our Lemma \ref{lem B2}.

\begin{lemma}\label{RS06-thm2}
  Let $\beta\in \mathbb{R}^{*}$, $\beta\neq 1$, be fixed.
  For   $N\ge 2$ and $X>0$, define
$$
\mathscr{B}_{0}
:=\Big| \Big\{(n_{1},n_{2},n_{3},n_{4})
\ \Big|\   n_{i} \sim N, i=1, 2, 3, 4, \
\Big|    \frac{n_{1}^{\beta}}{N^{\beta}} + \frac{n_{2}^{\beta}}{N^{\beta}}  -  \frac{n_{3}^{\beta}}{N^{\beta}}  - \frac{n_{4}^{\beta}}{N^{\beta}}  \Big|   \le\dfrac{1}{X}\Big\}    \Big|.
$$
Then
$$
\mathscr{B}_{0}\ll_{\epsilon}           N^{4+\epsilon}\bigg(\dfrac{1}{N^{2}}+\dfrac{1}{X}\bigg).
$$
\end{lemma}

Define
\begin{equation}  \label{uv}
\varphi_{n_{s},n_{t}}(m):=
\dfrac{N^{\gamma}}{n_{s}^{\gamma}+\mu(m)}-\dfrac{N^{\gamma}}{n_{t}^{\gamma}+\mu(m)},
\qquad
\psi_n( m):=\dfrac{N^{\gamma}}{n^{\gamma}+\nu(m)},
\end{equation}
and write
$$
 |\varphi_{n_{1},n_{2}}-\varphi_{n_{3},n_{4}}|:=\sup\limits_{  m_1, m_2\sim M}  | \varphi_{n_{1},n_{2}}(m_1)-\varphi_{n_{3},n_{4}}( m_2)|,
$$
$$
 |\psi_{n_1}-\psi_{n_2}|                      :=\sup\limits_{  m_1, m_2\sim M} | \psi_{n_1}(m_{1})-\psi_{n_1}(m_{2})|,
$$
where $0<\mu(m)<U\le 1$ and $0<\nu(m)<V\le 1$ for $ m\sim M$.

Applying Lemma \ref{RS06-thm2} we get the following lemma.

\begin{lemma}\label{lem B2}
	Let $\gamma>0$   be  fixed, $M, N, X\ge 1$ and $\varphi_{n_{s},n_{t}}$ be defined in \eqref{uv}.
Define
$$
\mathscr{B}_{2}:= | \{   (n_{1},n_{2},n_{3},n_{4})
\ |\
n_{i}\sim N,i=1,2,3,4,
 |\varphi_{n_{1},n_{2}} -\varphi_{n_{3},n_{4}} |  \le   {X}^{-1} \}  |.
$$
where $\varphi_{n_{s},n_{t}}$ is defined in \eqref{uv}.
If $X\le U^{-1}N^{\gamma}$,
then
$$
\mathscr{B}_{2}\ll_{\epsilon} N^{4+\epsilon}\bigg(\dfrac{1}{N^{2}}+\dfrac{1}{X}\bigg).
$$
\end{lemma}

\begin{proof}
 The condition in $\mathscr{B}_{2}$ is equivalent to
	\begin{equation}\label{sumcon2}
	\sup\limits_{m_{i}\sim M, n_i\sim N}\left|S_{1}Q_{1}-S_{2}Q_{2}\right|\le \dfrac{1}{X},
	\end{equation}
where
$$
S_{1}:=S(n_1, n_2):=\dfrac{N^{\gamma}}{n_{1}^{\gamma}}-\dfrac{N^{\gamma}}{n_{2}^{\gamma}},\qquad
Q_{1}:=Q(n_1, n_2, m_1):=\dfrac{n_{1}^{\gamma}n_{2}^{\gamma}}{(n_{1}^{\gamma}+\mu(m_1))(n_{2}^{\gamma}+\mu(m_{1}))},
$$
$$
S_{2}:=S(n_3, n_4):=\dfrac{N^{\gamma}}{n_{3}^{\gamma}}-\dfrac{N^{\gamma}}{n_{4}^{\gamma}}, \qquad
Q_{2}:=Q(n_3, n_4, m_2):=\dfrac{n_{3}^{\gamma}n_{4}^{\gamma}}{(n_{3}^{\gamma}+\mu(m_2))(n_{4}^{\gamma}+\mu(m_{2}))}
$$
temporarily in this proof.
It is easy to see that
\begin{equation}\label{3211}
 0<{1}/{4}\le        (1+U N^{-\gamma})^{-2}    \le    Q_{i}<1
\end{equation}
for  $i=1,2$ and
\begin{equation}\label{3222}
\sup\limits_{m_{1}\sim M}|S_{1}Q_{1}|
=\sup\limits_{m_{1} \sim M}
\dfrac{|n_{2}^{\gamma}-n_{1}^{\gamma}|N^{\gamma}}
{(n_{1}^{\gamma}+\mu(m_1))(n_{2}^{\gamma}+\mu(m_1))}
\le 1.
\end{equation}
Using the simple identity
$S_{1}-S_{2}= (S_{1}Q_{1}-S_{2}Q_{2})Q_{2}^{-1}+S_{1}Q_{1}(Q_{1}^{-1}-Q_{2}^{-1})$
and noticing that  $Q_{1},Q_{2}>0$,
we get
$$
\sup_{n_i\sim N}|S_{1}-S_{2}|\le \sup|S_{1}Q_{1}-S_{2}Q_{2}|(\inf Q_{2})^{-1}+\sup|S_{1}Q_{1}|\big((\inf Q_{1})^{-1}-(\sup Q_{2})^{-1}\big).
$$
Thus, if $n_{i}\sim N$, $i=1,2,3,4$, satisfy  the condition in $\mathscr{B}_{2}$,
then \eqref{sumcon2}, \eqref{3211}, \eqref{3222} and $X\le U^{-1}N^{\gamma}$ imply
$$
\sup\limits_{n_{i}\sim N}|S_{1}-S_{2}|
\le {X^{-1}}{(1+U  N^{-\gamma})^{2}}+    \big((1+U  N^{-\gamma})^{2}-1\big)
\le {7}{X}^{-1},
$$
which is a weaker condition than \eqref{sumcon2}.
	Finally, Lemma \ref{RS06-thm2}  completes the proof. %from asserts that we have the bound
%	\begin{equation}
%	\mathscr{B}_{2}\ll _{\epsilon} N^{4+\epsilon}\bigg(\dfrac{1}{N^{2}}+\dfrac{1}{X}\bigg).
%	\end{equation}	
\end{proof}

\begin{lemma}\label{lem B3}
Let $\gamma>0$   be  fixed, $M, N, X\ge 1$ and $\psi_{n}$ be defined in \eqref{uv}.
	Define
$$
\mathscr{B}_{3}
:= \big| \{(n_{1},n_{2})\
|\
n_{i}\sim N,i=1,2,
 |\psi_{n_1} -\psi_{n_2}   |     \le {X}^{-1} \} \big|,
$$
where $\psi_{n}$ is defined in \eqref{uv}.
If $X\le V^{-1}N^{\gamma}$,
then
$$
\mathscr{B}_{3}\ll N^{2}\bigg(\dfrac{1}{N}+\dfrac{1}{X}\bigg).
$$
\end{lemma}

\begin{proof}
For $i=1,2$, write	
$$
S_{i}:=S(n_i):=\dfrac{N^{\gamma}}{n_{i}^{\gamma}},\qquad
Q_{i}:=Q(n_i, m_i):=\dfrac{n_{i}^{\gamma}}{n_{i}^{\gamma}+\nu(m_i)}
$$
%$$
%S_{2}:=S(n_2):=\dfrac{N^{\gamma}}{n_{2}^{\gamma}}, \qquad
%Q_{2}:=Q(n_2, m_2):=\dfrac{n_{2}^{\gamma}}{n_{2}^{\gamma}+\nu(m_2)}
%$$
temporarily in this proof.
Then the   condition in $\mathscr{B}_{3}$ is equivalent to
	\begin{equation}\label{sumcon3}
	\sup\limits_{m_{i}\sim M, n_i\sim N}\left|S_{1}Q_{1}-S_{2}Q_{2}\right|\le  {X}^{-1}.
	\end{equation}
	An argument similar to the one used in the proof of Lemma \ref{lem B2} shows that
$|S_{1}-S_{2}|\le 3/X$, if \eqref{sumcon3} is satisfied.
Thus our problem reduces to the following diophantine problem
$$	%\begin{equation}\label{Bstar}
%	\mathscr{B}_{*}=
{\underset{|(n_{1}/N)^{-\gamma}-(n_{2}/N)^{-\gamma}|\le 3/X}{\sum\limits_{n_{1},n_{2}\sim N}}}1.
%\textcolor[rgb]{0.00,0.00,1.00}{=
%{\underset{|n_{1}-n_{2}|c(\gamma)N^{-1} \le 3/X}{\sum\limits_{n_{1},n_{2}\sim N}}}1.
%=
%{\underset{|n_{1}-n_{2}| \textcolor[rgb]{1.00,0.00,0.00}{\le cN/X}}{\sum\limits_{n_{1},n_{2}\sim N}}}1.}
$$%	\end{equation}
	Using Lagrange mean value theorem, we have 
$$
 2^{-\gamma-1}\gamma \frac{|n_1-n_2|}{N}   
 \le
 \Big|  \Big(  \frac{n_{1}}{N} \Big)^{-\gamma}    -    \Big(\frac{n_{2}}{N}\Big)^{-\gamma}\Big| 
\le 
\gamma  \frac{|n_1-n_2|}{N} .
$$
Thus, there are $O(1+N/X)$ eligible $n_{2}$ for each fixed $n_{1}$,
which complete the proof.
% have	
%	\begin{equation}
%	\mathscr{B}_{3}\ll N^{2}\bigg(\dfrac{1}{N}+\dfrac{1}{X}\bigg).
%	\end{equation}
\end{proof}

\subsection{Proof of Theorem \ref{thm_3DES}}
 	We write the exponential sum in \eqref{Sdelta} as
	\begin{equation}\label{Sd}
S_{\delta}:=	S_{\delta}(H,M,N)=\sum\limits_{h\sim H}\sum\limits_{m\sim M}a(h,m)\sum\limits_{n\sim N}b(n)e\left(X\dfrac{H^{-\alpha}M^{\beta}}{h^{-\alpha}m^{\beta}}\dfrac{N^{\gamma}}{n^{\gamma}+\delta m^{-\beta}}\right).
	\end{equation}
Firstly, suppose $X\gg HM$.	
By Cauchy's inequality, we have
	\begin{equation}\label{Sd2}
	\begin{aligned}
	|S_{\delta}|^{2}&\le HM\sum\limits_{h\sim H}\sum\limits_{m\sim M}\left|\sum\limits_{n\sim N}b(n)e\left(X\dfrac{H^{-\alpha}M^{\beta}}{h^{-\alpha}m^{\beta}}\dfrac{N^{\gamma}}{n^{\gamma}+\delta m^{-\beta}}\right)\right|^{2}
	\\&\le HM\sum\limits_{h\sim H}\sum\limits_{m\sim M}\sum\limits_{n_{1}\sim N}\sum\limits_{n_{2}\sim N}b(n_{1})\overline{b(n_{2})}e\bigg(X\dfrac{H^{-\alpha}M^{\beta}}{h^{-\alpha}m^{\beta}}\varphi_{n_{1},n_{2}}(m)\bigg),
	\end{aligned}
	\end{equation}
where   $\varphi_{n_{1},n_{2}}(m)$ is defined in \eqref{uv} with $\mu(m)=\delta m^{-\beta}.$
Set
$$
\mathscr{X}=\left\{\varphi_{n_{1},n_{2}}  \mid n_{1},n_{2}\sim N\right\},
\qquad
\mathscr{Y}=\left\{X\dfrac{H^{-\alpha}M^{\beta}}{h^{-\alpha}m^{\beta}}  \ \bigg|\ h\sim H,m\sim M\right\}.
$$
	It is easy to see that $\sup\mathscr{X}\le 1$ and $\sup\mathscr{Y}\le 2^\alpha X$.
Note that for $m_1, m_2\sim M$,
\begin{equation*}
	\begin{aligned}	
|\varphi_{n_{1},n_{2}}(m_1)-\varphi_{n_{1},n_{2}}(m_2)| &=\left|\dfrac{N^{\gamma}}{n_{1}^{\gamma}+\mu(m_{1})}-\dfrac{N^{\gamma}}{n_{2}^{\gamma}+\mu(m_{1})}-\dfrac{N^{\gamma}}{n_{1}^{\gamma}+\mu(m_{2})}+\dfrac{N^{\gamma}}{n_{2}^{\gamma}+\mu(m_{2})}\right|
	\\&\le \left|\dfrac{N^{\gamma}}{n_{1}^{\gamma}+\mu(m_{1})}-\dfrac{N^{\gamma}}{n_{1}^{\gamma}+\mu(m_{2})}\right|+\left|\dfrac{N^{\gamma}}{n_{2}^{\gamma}+\mu(m_{1})}-\dfrac{N^{\gamma}}{n_{2}^{\gamma}+\mu(m_{2})}\right|
	\\&\le\dfrac{N^{\gamma}|\mu(m_{2})-\mu(m_{1})|}{n_{1}^{2\gamma}}+\dfrac{N^{\gamma}|\mu(m_{2})-\mu(m_{1})|}{n_{2}^{2\gamma}}
	\\&\ < \dfrac{2\delta M^{-\beta}}{N^{\gamma}} \le \dfrac{ {K}}{4X},
	\end{aligned}
	\end{equation*}
since $X\le   (8\delta)^{-1}KM^\beta N^\gamma$.
So we can apply Proposition \ref{prop2.2} to \eqref{Sd2} to get
	\begin{equation}\label{sum}
	|S_{\delta}|^{2}\ll HM K^{\frac{1}{2}} X^{\frac{1}{2}}\mathscr{B}_{1}^{\frac{1}{2}}\mathscr{B}_{2}^{\frac{1}{2}},
	\end{equation}
where $\mathscr{B}_{1}$ and $\mathscr{B}_{2}$ are as in Lemma \ref{lem B1} and Lemma \ref{lem B2}.
Applying  Lemma \ref{lem B1} and Lemma \ref{lem B2} with $U=\delta M^{-\beta}$, we get
$$
\mathscr{B}_{1}\ll_{\epsilon} (HM)^{2+\epsilon} \dfrac{1}{HM},
\qquad
 \mathscr{B}_{2}\ll_{\epsilon} N^{4+\epsilon}\left(\dfrac{1}{N^{2}}+\dfrac{K}{X}\right),
$$
since $X\gg HM$.
	Inserting this into \eqref{sum}, we get
\begin{equation}\label{4.3}
	S_{\delta}(H,M,N)
\ll_{\epsilon}(HMN)^{1+\epsilon}
\bigg(\Big(\dfrac{KX}{HMN^{2}}\Big)^{\frac{1}{4}}  +  \Big(\dfrac{K^2}{HM} \Big)^{\frac{1}{4}}\bigg)%+\dfrac{1}{N^{\frac{1}{2}}}+\dfrac{1}{X^{\frac{1}{4}}},
	\end{equation}
	which implies \eqref{SdeltaR} for $X\gg HM$.
	
Now, we suppose $X\ll HM$.
We will directly apply   the generalized double large sieve inequality we get in Proposition \ref{prop2.2} to \eqref{Sd}.
Set
$$
\mathscr{X}=\{\psi_n  \mid   n\sim N\},
\qquad
\mathscr{Y}=\left\{X\dfrac{H^{-\alpha}M^{\beta}}{h^{-\alpha}m^{\beta}}
\ \bigg|\
h\sim H,m\sim M\right\},
$$
where   $\psi_n$ is defined in \eqref{uv} with $\nu(m)=\delta m^{-\beta}.$
It is easy to  see that $\sup\mathscr{X}\le 1$ and $\sup\mathscr{Y}\le 2^\alpha X$.
Note that for $m_1, m_2\sim M$,
\begin{equation*}
	\begin{aligned}	
|\psi_{n}(m_1)-\psi_{n}(m_2)|
%&=\left|\dfrac{N^{\gamma}}{n^{\gamma}+\mu(m_{1})} -\dfrac{N^{\gamma}}{n^{\gamma}+\mu(m_{2})} \right|
\le\dfrac{N^{\gamma}|\nu(m_{2})-\nu(m_{1})|}{n^{2\gamma}}
< \dfrac{2\delta M^{-\beta}}{N^{\gamma}}  \le \dfrac{K}{4X},
	\end{aligned}
	\end{equation*}
 since $X\le (8\delta)^{-1}KM^\beta N^\gamma$. %=o(M^{\beta}N^{\gamma})$.
Then \eqref{condition} is satisfied
so that we can apply  Proposition \ref{prop2.2} to \eqref{Sd} to get
	\begin{equation}
	S_{\delta}\ll   K^{\frac{1}{2}}   X^{\frac{1}{2}}\mathscr{B}_{1}^{\frac{1}{2}}\mathscr{B}_{3}^{\frac{1}{2}},
	\end{equation}
where  $\mathscr{B}_{1}$ and $\mathscr{B}_{3}$ are as in Lemma \ref{lem B1} and Lemma \ref{lem B3}.
Applying  Lemma \ref{lem B1} and Lemma \ref{lem B3} with $V=\delta M^{-\beta}$, we get
$$
\mathscr{B}_{1}\ll_{\epsilon} (HM)^{2+\epsilon}\dfrac{1}{X},
\qquad
\mathscr{B}_{3}\ll N^{2}\left(\dfrac{1}{N}+\dfrac{K}{X}\right),
$$
since $X\ll HM$.
Therefore,
	\begin{equation}\label{4.5}
	S_{\delta}(H,M,N)
	\ll_{\epsilon}(HMN)^{1+\epsilon}\bigg(     \Big(\dfrac{K}{N}\Big)^{\frac{1}{2}}+  \dfrac{K}{X^{\frac{1}{2}}}   \bigg).
	\end{equation}
Combining \eqref{4.3} and \eqref{4.5}, we get Theorem \ref{thm_3DES}.

\section{Proof of Theorem  \ref{thm_lambda}}

To prove Theorem  \ref{thm_lambda}, we need to cite some useful lemmas and then prove a key inequality in Proposition \ref{S_delta}.

\subsection{Preliminary lemmas}

In this subsection, we shall cite three lemmas, which will be needed in the the proof of Proposition \ref{S_delta}.
The first one is a direct corollary of \cite[Proposition 3.1]{LWY1}.

\begin{lemma}\label{LWY3DES}
	Let $\alpha,\beta,\gamma>0$ and $\delta\in\mathbb{R}$ be some constants. For $X>0$ and $H,M,N\ge 1$, define
	\begin{equation*}
	S_{\delta}
	=S_{\delta}(H,M,N):
	=\sum\limits_{h\sim H}\sum\limits_{m\sim M}\sum\limits_{n\sim N}a_{h,n}b_{m}e\bigg(X\dfrac{M^{\beta}N^{\gamma}}{H^{\alpha}}\dfrac{h^{\alpha}}{m^{\beta}n^{\gamma}+\delta}\bigg),
	\end{equation*}
	where $a_{h,m}, b_{n}\in \mathbb{C}$ such that $|a_{h,m}|,|b_{n}|\le 1$.
	For any $\epsilon>0$, we have
	\begin{equation}\label{Wu3DES}
	S_{\delta}\ll ((X^{\kappa}H^{2+\kappa}M^{2+\kappa}N^{1+\kappa+\lambda})^{1/(2+2\kappa)}+HMN^{1/2}+(HM)^{1/2}N+X^{-1/2}HMN)X^{\epsilon}
	\end{equation}
	uniformly for  $H\le M^{\beta-1}N^{\gamma}$ and $0\le \delta\le {\epsilon^{-1}}$,
where $(\kappa,\lambda)$ is an exponent pair and the implied constant depends on $(\alpha,\beta,\gamma,\epsilon)$ only.
\end{lemma}

The second one is due to Vaughan (\cite[(3)]{V80}).

\begin{lemma}\label{Vaughan}
	There are six real arithmetical functions $\alpha_{k}(n)$ verifying $|\alpha_{k}(n)|\ll_{\epsilon}n^{\epsilon}$
	for	$n>1$ $(1\le k\le 6)$ such that, for all $D>100$ and any arithmetical function $g$, we have
	\begin{equation}\label{Vaughan=}
	\sum\limits_{D<d\le 2D}\Lambda(d)g(d)=S_{1}+S_{2}+S_{3}+S_{4},
	\end{equation}
	where
\begin{align*}
& S_{1}:=\sum\limits_{m\le D^{1/3}}\alpha_{1}(m)\sum\limits_{D<mn\le 2D}g(mn),   \\
& S_{2}:=\sum\limits_{m\le D^{1/3}}\alpha_{2}(m)\sum\limits_{D<mn\le 2D}g(mn)\log n,    \\
& S_{3}:=\sum\limits_{\substack{D^{1/3}<m,n\le D^{2/3}\\D<mn\le 2D}}\alpha_{3}(m)\alpha_{4}(n)\sum\limits_{D<mn\le 2D}g(mn),  \\
& S_{4}:=\sum\limits_{\substack{D^{1/3}<m,n\le D^{2/3}\\D<mn\le 2D}}\alpha_{5}(m)\alpha_{6}(n)\sum\limits_{D<mn\le 2D}g(mn).
\end{align*}
The sums $S_{1}$ and $S_{2}$ are called as type I, $S_{3}$ and $S_{4}$ are called as type II.
\end{lemma}

The third one is due to Vaaler (\cite[Theorem A.6]{VDC}).
\begin{lemma}\label{psiVaaler}
	Let $\psi(t)=\{t\}-\frac{1}{2}$, where $\{t\}$ means the fractional part of the real number $t$.
	For $x\ge1$ and $H\ge1$, we have \begin{equation}\label{psi1}
	\psi(x)=-\sum\limits_{1\le |h|\le H}\Phi\bigg(\dfrac{h}{H+1}\bigg)\dfrac{e(hx)}{2\pi ih}+R_{H}(x),
	\end{equation}
	where  $\Phi(t):=\pi t(1-|t|)\cot (\pi t)+|t|$, and the error term $R_{H}(x)$ satisfies
	\begin{equation}\label{psi2}
	|R_{H}(x)|\le \dfrac{1}{2H+2}\sum\limits_{0\le|h|\le H}\bigg(1-\dfrac{|h|}{H+1}\bigg)e(hx).
	\end{equation}
\end{lemma}

\subsection{Two key inequalities}
Denote
\begin{equation}
	\mathfrak{S}_{\delta}(x,D):=\sum\limits_{D<d\le 2D}\Lambda(d)\psi\bigg(\dfrac{x}{d+\delta}\bigg),
	\end{equation}
where $\delta$ is a non-negative constant.	
For $\mathfrak{S}_{\delta}(x,D)$, \cite[Proposition 4.1]{LWY1} gives the following bound
which will be used in the proof of Theorem \ref{thm_lambda} in the case that $D$ is small.

\begin{lemma}
		Let $\delta\ge 0$ be a  constant.
		Then
		\begin{equation}\label{WuBound}
		\mathfrak{S}_{\delta}(x,D)\ll_{\epsilon}(x^{2}D^{7})^{1/12}x^{\epsilon}
		\end{equation}
		uniformly for $x\ge 3$ and $x^{6/13}\le D\le x^{2/3}$.
\end{lemma}

We still need another bound  of $\mathfrak{S}_{\delta}(x,D)$ in the proof of Theorem \ref{thm_lambda} in the case that $D$ is  big.

	\begin{proposition}\label{S_delta}
		Let $\delta\ge 0$ be a  constant.
		Then
		\begin{equation}\label{LMbound}
		\mathfrak{S}_{\delta}(x,D)\ll_{\epsilon}(  D^{17/19}  +  x^{1/6}D^{329/570}  )x^{\epsilon}
		\end{equation}
		uniformly for $x\ge 3$ and   $x^{11/21}\le D\le x^{3/4}$.
	\end{proposition}

	\begin{proof}
		We apply the Vaughan identity \eqref{Vaughan=} with $g(d)=\psi\bigg(\dfrac{x}{d+\delta}\bigg)$ to write
			\begin{equation}\label{SSSS}
			\mathfrak{S}_{\delta}(x,D)=\mathfrak{S}_{\delta,1}+\mathfrak{S}_{\delta,2}+\mathfrak{S}_{\delta,3}+\mathfrak{S}_{\delta,4},
			\end{equation}
			where
\begin{align*}
& \mathfrak{S}_{\delta,1}:=\sum\limits_{m\le D^{1/3}}\alpha_{1}(m)\sum\limits_{D<mn\le 2D}\psi\bigg(\dfrac{x}{mn+\delta}\bigg),     \\
&  \mathfrak{S}_{\delta,2}:=\sum\limits_{m\le D^{1/3}}\alpha_{2}(m)\sum\limits_{D<mn\le 2D}\psi\bigg(\dfrac{x}{mn+\delta}\bigg)\log n,\\
&  \mathfrak{S}_{\delta,3}:=\sum\limits_{\substack{D^{1/3}<m,n\le D^{2/3}\\D<mn\le 2D}}\alpha_{3}(m)\alpha_{4}(n)\sum\limits_{D<mn\le 2D}\psi\bigg(\dfrac{x}{mn+\delta}\bigg),\\
& \mathfrak{S}_{\delta,4}:=\sum\limits_{\substack{D^{1/3}<m,n\le D^{2/3}\\D<mn\le 2D}}\alpha_{5}(m)\alpha_{6}(n)\sum\limits_{D<mn\le 2D}\psi\bigg(\dfrac{x}{mn+\delta}\bigg).
\end{align*}

A. Estimation of $\mathfrak{S}_{\delta,3}$ and $\mathfrak{S}_{\delta,4}$.

		Using Lemma \ref{psiVaaler} and splitting the interval of summation into dyadic intervals, we can
		write
		\begin{equation}\label{Sdelta3}
		\begin{split}
		\mathfrak{S}_{\delta,3}
		=&-\dfrac{1}{2\pi i}
		\sum\limits_{H'}\sum\limits_{M}\sum\limits_{N}(\mathfrak{S}^\flat_{\delta,3}(H',M,N)+\overline{\mathfrak{S}^\flat_{\delta,3}(H',M,N)}) \\
		& +\sum\limits_{M}\sum\limits_{N}\mathfrak{S}^\dag_{\delta,3}(M,N)
		\end{split}
		\end{equation}
with
		$$
		\mathfrak{S}^\flat_{\delta,3}(H',M,N)
:=\dfrac{1}{H'}\sum\limits_{h\sim H'}a_{h}\underset{D<mn\le 2D}{\sum\limits_{m\sim M}\sum\limits_{n\sim N}}\alpha_{3}(m)\alpha_{4}(n)e\bigg(\dfrac{hx}{mn+\delta}\bigg),
		$$
		$$
		\mathfrak{S}^\dag_{\delta,3}(M,N)
:=\underset{D<mn\le 2D}{\sum\limits_{m\sim M}\sum\limits_{n\sim N}}\alpha_{3}(m)\alpha_{4}(n)R_{H}\bigg(\dfrac{x}{mn+\delta}\bigg),
		$$
		where $1\le H'\le H\le D^{2/19}$, $a_{h}:=\frac{H'}{h}\Phi (\frac{h}{H+1})\ll 1$, $MN\asymp D$,
		and  $D^{1/3}\le M\le D^{1/2}\le N\le D^{2/3}$ in view of the symmetry of the variables $m$ and $n$.

		Firstly, we bound $\mathfrak{S}^\flat_{\delta,3}(H',M,N)$.
We can remove the  multiplicative condition $D<mn\le 2D$ at the cost of a factor $\log x$ to get
\begin{equation}\label{S_delta3 f}
\mathfrak{S}^\flat_{\delta,3}(H',M,N)
		\ll \dfrac{x^{\epsilon}}{H'}\sum\limits_{h\sim H'}\sum\limits_{m\sim M}\sum\limits_{n\sim N} a_{h}\dfrac{\alpha_{3}(m)}{M^{\epsilon}}\dfrac{\alpha_{4}(n)}{N^{\epsilon}}e\bigg(\dfrac{hx}{mn+\delta}\bigg).
\end{equation}
Let $K=D^{7/380}$ and $1/3<\theta<1/2$ be a parameter to be chosen later.
For the case $D^{\theta}\le M\le D^{1/2}\le N\le D^{1-\theta}$,
note that $H\le K D^{13/380}$ and $x^{11/21}\le D$ imply that $xH'/MN%\ll xH/D\ll xKD^{-367/380}
=o(DK)$,
then  %$X=o(KMN)$ and
we can  apply Theorem \ref{thm_3DES} with
$\alpha=\beta=\gamma=1$, $K=D^{7/380}$,
$(X,H,M,N)=({xH'}/{MN},H',M,N)$  to get
		\begin{align}\label{Sdeltaflat}
		\mathfrak{S}^\flat_{\delta,3}(H',M,N) \nonumber
%		&\ll \dfrac{x^{\epsilon}}{H'}\sum\limits_{h\sim H'}\sum\limits_{m\sim M}\sum\limits_{n\sim N} a_{h}\dfrac{\alpha_{3}(m)}{M^{\epsilon}}\dfrac{\alpha_{4}(n)}{N^{\epsilon}}e\bigg(\dfrac{hx}{mn+\delta}\bigg)\nonumber
		& \ll \big(    x^{\frac{1}{4}}K^{\frac{1}{4}}M^{\frac{1}{2}}N^{\frac{1}{4}}   +H'^{-\frac{1}{4}}K^{\frac{1}{2}}M^{\frac{3}{4}}N   +  K^{\frac{1}{2}}MN^{\frac{1}{2}}
		                +x^{-\frac{1}{2}} K  H'^{-\frac{1}{2}}M^{\frac{3}{2}}N^{\frac{3}{2}}\big)x^{\epsilon}
		\\
&\ll \big(    x^{1/4}D^{3/8} K^{{1}/{4}}   +  D^{1-\theta/4} K^{{1}/{2}}  +  D^{ {3}/{4}}K^{{1}/{2}}
		      +x^{-{1}/{2}}D^{3/2}K\big)x^{\epsilon},
		\end{align}
where we have used $H'\ge 1$, $MN\asymp D$ and $D^{\theta}\le M\le D^{1/2}\le N\le D^{1-\theta}$.
For the other case $D^{1/3}\le M\le D^{\theta}\le D^{1-\theta}\le N\le D^{2/3}$,
apply Lemma \ref{LWY3DES} with $\alpha=\beta=\gamma=1$,
$(X,H,M,N)=(xH'/MN,H',M,N)$ and $(\kappa, \lambda)=({1}/{2}, {1}/{2})$ to \eqref{S_delta3 f} to get
\begin{equation}
\begin{split}\label{Sdeltaflat2}
\mathfrak{S}^\flat_{\delta,3}(H',M,N)
		%&\ll \dfrac{x^{\epsilon}}{H'}\sum\limits_{h\sim H'}\sum\limits_{m\sim M}\sum\limits_{n\sim N} a_{h}\dfrac{\alpha_{3}(m)}{M^{\epsilon}}\dfrac{\alpha_{4}(n)}{N^{\epsilon}}e\bigg(\dfrac{hx}{mn+\delta}\bigg)\nonumber\\
& \ll \big(x^{1/6}M^{2/3}N^{1/2}     +  MN^{ {1}/{2}}    +  M^{1/2}N   +  		x^{-{1}/{2}}D^{ {{3}/{2}}}  \big)x^{\epsilon}
		\\&\ll \big(   x^{1/6}D^{1/2+\theta/6}   +   D^{ 5/6}   \big)x^{\epsilon},
		\end{split}
\end{equation}
where we have used $M\le D^{\theta}$, $N\le D^{2/3}$, $H'\ge 1$, $D\le x^{ {3/4}}$.
Summarize \eqref{Sdeltaflat} and \eqref{Sdeltaflat2}, we get
\begin{equation}\label{S_delta3 flat}
\mathfrak{S}^\flat_{\delta,3}(H',M,N)
\ll \big(
        x^{\frac{1}{4}}D^{\frac{3}{8}} K^{\frac{1}{4}}   +  D^{1-\frac{\theta}{4}}K^{\frac{1}{2}}
		 +  x^{\frac{1}{6}}D^{\frac{1}{2}  + \frac{\theta}{6}}   +D^{\frac{ 8}{9}}
        \big)x^{\epsilon},
\end{equation}		
where we have used $1\le x^{11/21}\le D\le x^{3/4}$.

Secondly, we bound $\mathfrak{S}^\dag_{\delta,3}(M,N)$.
Using \eqref{psi2}, we have
\begin{equation*}
		\begin{aligned}
\mathfrak{S}^\dag_{\delta,3}(M,N)
%		& \ll \dfrac{x^{\epsilon}}{H}\sum\limits_{m\sim M}\sum\limits_{n\sim N}\sum\limits_{0\le |h|\le H}\bigg(1-\dfrac{|h|}{H+1}\bigg)e\bigg(\dfrac{hx}{mn+\delta}\bigg)\\
		& \ll \dfrac{x^{\epsilon}}{H}\bigg(MN+\sum\limits_{1\le |h|\le H}\sum\limits_{m\sim M}\sum\limits_{n\sim N} \bigg(1-\dfrac{|h|}{H+1}\bigg)e\bigg(\dfrac{hx}{mn+\delta}\bigg)\bigg)
\\
& \ll \dfrac{x^{\epsilon}}{H}\Big(MN+\max\limits_{1\le H'\le H}\big|\tilde{\mathfrak{S}}^{\dag}_{\delta,3}(H',M,N)\big|\Big),
		\end{aligned}
\end{equation*}
		where
$$
		\tilde{\mathfrak{S}}^{\dag}_{\delta,3}(H',M,N):
		=\sum\limits_{h\sim H'}\sum\limits_{m\sim M}\sum\limits_{n\sim N}\bigg(1-\dfrac{|h|}{H+1}\bigg)e\bigg(\dfrac{hx}{mn+\delta}\bigg).
$$
Clearly we can bound $\tilde{\mathfrak{S}}^{\dag}_{\delta,3}(H',M,N)$ in the same way as $\mathfrak{S}^\flat_{\delta,3}(H',M,N)$ and  get
\begin{equation}\label{Sdag}
\begin{split}
&\mathfrak{S}^\dag_{\delta,3}(M,N)
\ll \big(
        DH^{-1}  +
        x^{\frac{1}{4}}D^{\frac{3}{8}} K^{\frac{1}{4}}   +  D^{1-\frac{\theta}{4}}K^{\frac{1}{2}}
	        + x^{\frac{1}{6}}D^{\frac{1}{2}+\frac{\theta}{6}}   +D^{\frac{8}{9}}
         \big)x^{\epsilon}.
\end{split}
\end{equation}

Inserting \eqref{S_delta3 flat} and \eqref{Sdag} into \eqref{Sdelta3} gives
\begin{equation*}%\label{Sdeltadag}
\mathfrak{S}_{\delta,3}
\ll \big(
DH^{-1}  +
        x^{\frac{1}{4}}D^{\frac{3}{8}} K^{\frac{1}{4}}   +  D^{1-\frac{\theta}{4}}K^{\frac{1}{2}}
	        + x^{\frac{1}{6}}D^{\frac{1}{2}+\frac{\theta}{6}}   +D^{\frac{8}{9}}
\big)x^{\epsilon}.
\end{equation*}
Taking $H=D^{2/19}$,  $K=D^{7/380}$  and $\theta=44/95$ gives
\begin{equation}\label{S_delta3}
\mathfrak{S}_{\delta,3}
\ll \big(D^{17/19}  +  x^{1/6}D^{329/570}\big)x^{\epsilon},
\end{equation}
where we have used $1\le x^{11/21}\le D\le x^{3/4}$.

Clearly the same estimate also holds for $\mathfrak{S}_{\delta,4}$.
		
B. Estimation of $\mathfrak{S}_{\delta,1}$ and $\mathfrak{S}_{\delta,2}$.

Since a simple partial summation allows us to treat %remove the smooth factor $\log n$ in $\mathfrak{S}_{\delta, 2}$,
$\mathfrak{S}_{\delta, 1}$ and $\mathfrak{S}_{\delta, 2}$   in the same way.
So,  we shall bound only $\mathfrak{S}_{\delta, 1}$.
Lemma \ref{psiVaaler} allows us to write
\begin{equation}\label{proof:prop_11}
\mathfrak{S}_{\delta, 1}
= - \frac{1}{2\pi\text{i}} \big(\mathfrak{S}^{\flat}_{\delta, 1} + \overline{\mathfrak{S}^{\flat}_{\delta, 1}}\big)
+ \mathfrak{S}^{\dagger}_{\delta, 1},
\end{equation}
where
\begin{align*}
\mathfrak{S}^{\flat}_{\delta, 1}
& := \sum_{m\le D^{1/3}} \alpha_1(m) \sum_{D<mn\le 2D}
\sum_{1\le l\le L}\frac{1}{l} \Phi \Big(\frac{l}{L+1} \Big) \text{e}\Big(\frac{lx}{mn+\delta} \Big),
\\
\mathfrak{S}^{\dagger}_{\delta, 1}
& := \sum_{m\le D^{1/3}} \alpha_1(m) \sum_{D<mn\le 2D} R_L\Big(\frac{x}{mn+\delta}\Big),
\end{align*}
and $L\ge 1$.
Applying the exponent pair $(\kappa',\lambda')$ to the sum over $n$,
we can derive
\begin{equation}\label{proof:prop_12}
\begin{aligned}
\mathfrak{S}^{\flat}_{\delta, 1}
& \ll x^{\varepsilon} \sum_{m\le D^{1/3}} \sum_{1\le l\le L} \frac{1}{l}
\bigg(  \Big(\frac{lx}{(D^2/m)}\Big)^{\kappa'} \Big(\frac{D}{m}\Big)^{\lambda'} + \frac{(D^2/m)}{lx}  \bigg)
\\
& \ll x^{\varepsilon} \sum_{m\le D^{1/3}}
\big(x^{\kappa'} D^{-2\kappa'+\lambda'} m^{\kappa'-\lambda'} L^{\kappa'} + x^{-1} D^2m^{-1}\big)
\\
& \ll \big(x^{\kappa'} D^{(-5\kappa'+2\lambda'+1)/3} L^{\kappa'} + x^{-1} D^2\big) x^{\varepsilon}
\end{aligned}
\end{equation}
for all $L\ge 1$.
To bound $\mathfrak{S}^{\dagger}_{\delta, 1}$, apply \eqref{psi2} we get
\begin{align*}
\big|\mathfrak{S}^{\dagger}_{\delta, 1}\big|
& \ll x^{\varepsilon} \sum_{m\le D^{1/3}} \sum_{D<mn\le 2D}
\Big|R_L\Big(\frac{x}{mn+\delta}\Big)\Big|
\ll \big(DL^{-1} + \big|\widetilde{\mathfrak{S}}^{\dagger}_{\delta, 1}\big|\big) x^{\varepsilon},
\end{align*}
where
$$
\widetilde{\mathfrak{S}}^{\dagger}_{\delta, 1}
:= \frac{1}{L} \sum_{m\le D^{1/3}} \sum_{D<mn\le 2D}
\sum_{1\le |l|\le L} \Big(1-\frac{|l|}{L+1}\Big) {\rm e}\Big(\frac{lx}{mn+\delta}\Big).
$$
Clearly we can treat $\widetilde{\mathfrak{S}}^{\dagger}_{\delta, 1}$
in the same way as $\mathfrak{S}^{\flat}_{\delta, 1}$
and obtain
\begin{equation}\label{proof:prop_13}
\mathfrak{S}^{\dagger}_{\delta, 1}
\ll \big(D L^{-1}  +  x^{\kappa'} D^{(-5\kappa'+2\lambda'+1)/3} L^{\kappa'} + x^{-1} D^2\big) x^{\varepsilon}
\end{equation}
for $L\ge 1$.
Inserting \eqref{proof:prop_12} and \eqref{proof:prop_13} into \eqref{proof:prop_11}, it follows that
$$
\mathfrak{S}_{\delta, 1}
\ll \big(D L^{-1} + x^{\kappa'} D^{(-5\kappa'+2\lambda'+1)/3} L^{\kappa'} + x^{-1} D^2\big) x^{\varepsilon}
$$
for $L\ge 1$.
Optimising $L$ over $[1, \infty)$, we get
\begin{equation*}%\label{S_delta1}
\mathfrak{S}_{\delta, 1}
\ll \big((x^{3\kappa'} D^{-2\kappa'+2\lambda'+1})^{1/(3\kappa'+3)} + x^{\kappa'} D^{(-5\kappa'+2\lambda'+1)/3} + x^{-1} D^2\big) x^{\varepsilon}
\end{equation*}
The same estimate   also holds for $\mathfrak{S}_{\delta, 2}$.
Taking $ (\kappa', \lambda')=(\frac{1}{2}, \frac{1}{2})$ implies
\begin{equation}\label{S_delta12}
		\mathfrak{S}_{\delta,j}\ll \big(x^{1/3}D^{2/9}+x^{1/2}D^{-1/6}+x^{-1}D^{2}\big)x^{\epsilon}
\ll D^{17/19}
\end{equation}
for $j=1,2$, since $1\le x^{11/21}\le D\le x^{3/4}$.

Inserting \eqref{S_delta12} and \eqref{S_delta3} into \eqref{SSSS} completes the proof.  
\end{proof}

\subsection{Proof of Theorem \ref{thm_lambda}}
	
Let $E\in[x^{8/17},x^{1/2}]$ be a parameter to be chosen later.
	Write
	\begin{equation}\label{ThmS}
	\sum\limits_{n\le x}\Lambda\bigg(\bigg[\dfrac{x}{n}\bigg]\bigg):=S_{1}(x)+S_{2}(x)
	\end{equation}
	with
	$$
S_{1}(x)
	:=\sum\limits_{n\le E}\Lambda\bigg(\bigg[\dfrac{x}{n}\bigg]\bigg),\qquad
S_{2}(x)
    :=\sum\limits_{E< n\le x}\Lambda\bigg(\bigg[\dfrac{x}{n}\bigg]\bigg).
	$$
	Noticing that $\Lambda(n)\ll_{\epsilon} n^{\epsilon}$, we get
	\begin{equation}
	S_{1}(x)\ll_{\epsilon} Ex^{\epsilon}.
	\end{equation}
	
	Next we bound $S_{2}(x)$. Putting $d=[x/n]$, then $\frac{x}{n}-1<d\le \frac{x}{n}\Leftrightarrow \frac{x}{d+1}<n\le \frac{x}{d}$, thus we can write
	
	\begin{equation*}
	\begin{aligned}
	S_{2}(x)&=\sum\limits_{d\le x/E}\Lambda(d)\sum\limits_{x/(d+1)<n\le x/d}1\\
	%&=\sum\limits_{d\le x/N_{f}(k)}\Lambda(d)\bigg(\dfrac{x}{d}-\bigg\{\dfrac{x}{d}\bigg\}-\dfrac{x}{d+1}+\bigg\{\dfrac{x}{d+1}\bigg\}\bigg)\\
	&=\sum\limits_{d\le x/E}\Lambda(d)\bigg(\dfrac{x}{d}-\psi\bigg(\dfrac{x}{d}\bigg)-\dfrac{x}{d+1}+\psi\bigg(\dfrac{x}{d+1}\bigg)\bigg)\\
	&=x\sum\limits_{d\ge 1}\dfrac{\Lambda(d)}{d(d+1)}  +  R_{1}(x,k)  -  R_{0}(x,k)  +  O_{\epsilon}(Ex^{\epsilon}),
	\end{aligned}
	\end{equation*}
	where we have used the following bounds
$$
x\sum\limits_{d>x/E}\dfrac{\Lambda(d)}{d(d+1)}\ll_{\epsilon}E x^{\epsilon},
\qquad
\sum\limits_{d\le E}\Lambda(d)\bigg(\psi\bigg(\dfrac{x}{d+1}\bigg)-\psi\bigg(\dfrac{x}{d}\bigg)\bigg)\ll_{\epsilon}Ex^{\epsilon},
$$
	and defined
$$
R_{\delta}(x):=\sum\limits_{E< d\le x/E}\Lambda(d)\psi\bigg(\dfrac{x}{d+\delta}\bigg).
$$
Writing $D_{j}=x/(2^{j}E)$, we have $E< D_{j}\le x/E$ for $0\le j\le \frac{\log (x/E^{2})}{\log 2}$. Thus we can apply \eqref{WuBound} and \eqref{LMbound} to get
\begin{equation}\label{ThmRdelta}
\begin{aligned}
R_{\delta}(x)&=\sum\limits_{0\le j\le\log (x/E^2)/\log 2}\big|\mathfrak{S}_{\delta}(x,D_{j})\big|\\
&\ll_{\epsilon} \dfrac{\log (2x/E^{2})}{\log 2}\max\limits_{0\le j\le\log (x/E^2)/\log 2} \big|\mathfrak{S}_{\delta}(x,D_{j})\big|\\
&\ll_{\epsilon} \max\limits_{E<D\le x^{11/21}} \big|\mathfrak{S}_{\delta}(x,D)\big|x^{\epsilon}+ \max\limits_{x^{11/21}<D\le x/E} \big|\mathfrak{S}_{\delta}(x,D)\big|x^{\epsilon}
\\
&\ll_{\epsilon} \max\limits_{E<D\le x^{11/21}}    (x^2 D^7)^{1/12} x^{\epsilon}+ \max\limits_{x^{11/21}<D\le x/E} ( D^{17/19}  +  x^{1/6}D^{329/570}  )x^{\epsilon}
\\
&\ll_{\epsilon} (x^{17/36} +  x^{17/19}E^{-17/19}    +    x^{212/285}E^{-329/570})x^{\epsilon}.
\end{aligned}
\end{equation}
Take $E=x^{17/36}$, the proof is completed.

\vskip 5mm

\noindent{\bf Acknowledgements}.
This work is supported in part by the National Natural Science Foundation of China (Grant No. 11771252).

\end{document}